\documentclass[a4paper,9pt]{article}
\usepackage[cp1251]{inputenc}
\usepackage[english,russian]{babel}
\usepackage{amsmath}
\usepackage{amssymb}
\usepackage{amsfonts}
\def\udcs{512.5} 
\numberwithin{equation}{section}
\begin{document}
УДК \udcs \thispagestyle{empty}

\begin{center}
\textbf{Double-Layer Potentials for a Generalized Bi-Axially
Symmetric Helmholtz Equation II}
\end{center}

\begin{center}
\textit{Abdumauvlen Berdyshev$^{{1}}$, Anvar Hasanov$^{{2}}$and
Tuhtasin Ergashev$^{{3}}$ }

\end{center}

\bigskip
$^{{1}}$ Abai Kazakh National Pedagogical University, Almati,
Kazakhstan, \\
 $^{{2,3}}$Institute of Mathematics, Uzbek Academy of
Sciences, Tashkent, Uzbekistan

\bigskip

\textbf{Abstract:} The double-layer potential plays an important
role in solving boundary value problems for elliptic equations.
All the fundamental solutions of the generalized bi-axially
symmetric Helmholtz equation were known (\textit{Complex Variables
and Elliptic Equations}, \textbf{52}(8), 2007, 673--683.), and
only for the first one was constructed the theory of potential
(\textit{Sohag Journal of Mathematics} 2, No. 1, 1-10, 2015).
Here, in this paper, we aim at constructing theory of double-layer
potentials corresponding to the next fundamental solution. By
using some properties of one of Appell's hypergeometric functions
in two variables, we prove limiting theorems and derive integral
equations concerning a denseness of double-layer potentials.


\bigskip

\section{Introduction}

Potential theory has played a paramount role in both analysis and
computation for boundary value problems for elliptic partial
differential equations. Numerous applications can be found in
fracture mechanics, fluid mechanics, elastodynamics,
electromagnetic, and acoustics. Results from potential theory
allow us to represent boundary value problems in integral equation
form. For problems with known Green's functions, an integral
equation formulation leads to powerful numerical approximation
schemes.

The double-layer potential plays an important role in solving
boundary value problems of elliptic equations. The representation
of the solution of the (first) boundary value problem is sought as
a double-layer potential with unknown density and an application
of certain property leads to a Fredholm equation of the second
kind for determining the function (see [18] and [29]).

By applying a method of complex analysis (based upon analytic
functions), Gilbert [16] constructed an integral representation of
solutions of the following generalized bi-axially symmetric
Helmholtz equation:

\[
H_{\alpha ,\beta} ^{\lambda}  \left( {u} \right) \equiv u_{xx} +
u_{yy} + {\frac{{2\alpha}} {{x}}}u_{x} + {\frac{{2\beta}}
{{y}}}u_{y} - \lambda ^{2}u = 0, \quad \left( {0 < \alpha <
{\frac{{1}}{{2}}};0 < \beta < {\frac{{1}}{{2}}}} \right), \quad
\left( {H_{\alpha ,\beta} ^{\lambda}}   \right)
\]

Fundamental solutions of the equation $\left( {H_{\alpha ,\beta}
^{\lambda} } \right)$ were constructed recently (see [19]). In
fact, the fundamental solutions of the equation $\left( {H_{\alpha
,\beta} ^{\lambda}}   \right)$ when $\lambda = 0$ can be expressed
in terms of Appell's hypergeometric function in two variables of
the second kind, that is, the Appell function

\[
F_{2} \left( {a,b_{1} ,b_{2} ;c_{1} ,c_{2} ;x,y} \right)
\]

\noindent defined by (see [8, p.224,Eq.5.7.1(7)]; see also [4,
p.14,Eq.(12)] and [34, p.23, Eq.1.3(3)])

\begin{equation}
\label{eq1} F_{2} \left( {a;b_{1} ,b_{2} ;c_{1} ,c_{2} ;x,y}
\right) = {\sum\limits_{m,n = 0}^{\infty}  {}} {\frac{{\left( {a}
\right)_{m + n} \left( {b_{1}} \right)_{m} \left( {b_{2}}
\right)_{n}}} {{\left( {c_{1}}  \right)_{m} \left( {c_{2}}
\right)_{n} m!n!}}}x^{m}y^{n},
\end{equation}

\noindent where $\left( {a} \right)_{\nu}  $ denotes the general
Pochhammer symbol which is defined (for $a,\,\nu \in C\backslash
\{0\})$, in terms of the familiar Gamma function, by

\[
\left( {a} \right)_{\nu}  : = {\frac{{\Gamma \left( {a + \nu}
\right)}}{{\Gamma \left( {a} \right)}}} = {\left\{
{{\begin{array}{*{20}c}
 {1} \hfill & {(\nu = 0;\,\,\,a \in C\backslash \{0\})} \hfill \\
 {a(a + 1)...(a + \nu - 1)} \hfill & {(\nu = n \in N;\,\,\,a \in C\backslash
\{0\}),} \hfill \\
\end{array}}}  \right.}
\]

\noindent it being understood conventionally that $(0)_{0} : = 1$
and assumed tacitly and the $\Gamma - $quotient exists.

In case of$\,\lambda = 0$ fundamental solutions look like

\begin{equation}
\label{eq2} q_{1} \left( {x,y;x_{0} ,y_{0}}  \right) = k_{1}
\left( {r^{2}} \right)^{ - \alpha - \beta} F_{2} \left( {\alpha +
\beta ;\alpha ,\beta ;2\alpha ,2\beta ;\xi ,\eta}  \right),
\end{equation}

\begin{equation}
\label{eq3} q_{2} \left( {x,y;x_{0} ,y_{0}}  \right) = k_{2}
\left( {r^{2}} \right)^{\alpha - \beta - 1}x^{1 - 2\alpha}
x_{0}^{1 - 2\alpha}  F_{2} \left( {1 - \alpha + \beta ;1 - \alpha
,\beta ;2 - 2\alpha ,2\beta ;\xi ,\eta}  \right),
\end{equation}

\begin{equation}
\label{eq4} q_{3} \left( {x,y;x_{0} ,y_{0}}  \right) = k_{3}
\left( {r^{2}} \right)^{ - \alpha + \beta - 1}y^{1 - 2\beta}
y_{0}^{1 - 2\beta}  F_{2} \left( {1 + \alpha - \beta ;\alpha ,1 -
\beta ;2\alpha ,2 - 2\beta ;\xi ,\eta}  \right)
\end{equation}

\noindent and

\begin{equation}
\label{eq5} q_{4} \left( {x,y;x_{0} ,y_{0}}  \right) = k_{4}
\left( {r^{2}} \right)^{\alpha + \beta - 2}x^{1 - 2\alpha} y^{1 -
2\beta} x_{0}^{1 - 2\alpha}  y_{0}^{1 - 2\beta}  F_{2} \left( {2 -
\alpha - \beta ;1 - \alpha ,1 - \beta ;2 - 2\alpha ,2 - 2\beta
;\xi ,\eta}  \right),
\end{equation}

where

\begin{equation}
\label{eq6} k_{1} = {\frac{{2^{2\alpha + 2\beta}} }{{4\pi}}
}{\frac{{\Gamma \left( {\alpha}  \right)\Gamma \left( {\beta}
\right)\Gamma \left( {\alpha + \beta } \right)}}{{\Gamma \left(
{2\alpha}  \right)\Gamma \left( {2\beta} \right)}}},
\end{equation}

\begin{equation}
\label{eq7} k_{2} = {\frac{{2^{2 - 2\alpha + 2\beta}} }{{4\pi}}
}{\frac{{\Gamma \left( {1 - \alpha}  \right)\Gamma \left( {\beta}
\right)\Gamma \left( {1 - \alpha + \beta}  \right)}}{{\Gamma
\left( {2 - 2\alpha}  \right)\Gamma \left( {2\beta}  \right)}}},
\end{equation}

\begin{equation}
\label{eq8} k_{3} = {\frac{{2^{2 + 2\alpha - 2\beta}} }{{4\pi}}
}{\frac{{\Gamma \left( {\alpha}  \right)\Gamma \left( {1 - \beta}
\right)\Gamma \left( {1 + \alpha - \beta}  \right)}}{{\Gamma
\left( {2\alpha}  \right)\Gamma \left( {2 - 2\beta}  \right)}}},
\end{equation}

\begin{equation}
\label{eq9} k_{4} = {\frac{{2^{4 - 2\alpha - 2\beta}} }{{4\pi}}
}{\frac{{\Gamma \left( {1 - \alpha}  \right)\Gamma \left( {1 -
\beta}  \right)\Gamma \left( {2 - \alpha - \beta}
\right)}}{{\Gamma \left( {2 - 2\alpha}  \right)\Gamma \left( {2 -
2\beta}  \right)}}},
\end{equation}

\begin{equation}
\label{eq10} {\left. {{\begin{array}{*{20}c}
 {r^{2}} \hfill \\
 {r_{1}^{2}}  \hfill \\
 {r_{2}^{2}}  \hfill \\
\end{array}}}  \right\}}  = \left( {x{\begin{array}{*{20}c}
 { -}  \hfill \\
 { +}  \hfill \\
 { -}  \hfill \\
\end{array}} x_{0}}  \right)^{2} + \left( {y{\begin{array}{*{20}c}
 { -}  \hfill \\
 { -}  \hfill \\
 { +}  \hfill \\
\end{array}} y_{0} \,} \right)^{2},
\xi = {\frac{{r^{2} - r_{1}^{2}}} {{r^{2}}}}, \eta = {\frac{{r^{2}
- r_{2}^{2}}} {{r^{2}}}},
\end{equation}

The fundamental solutions (\ref{eq2}) - (\ref{eq5}) possess the
following properties:

\begin{equation}
\label{eq11} {\left. {x^{2\alpha} {\frac{{\partial q_{1} \left(
{x,y;x_{0} ,y_{0}} \right)}}{{\partial x}}}} \right|}_{x = 0} = 0,
\quad {\left. {y^{2\beta} {\frac{{\partial q_{1} \left( {x,y;x_{0}
,y_{0}} \right)}}{{\partial y}}}} \right|}_{y = 0} = 0,
\end{equation}

\begin{equation}
\label{eq12} {\left. {q_{2} \left( {x,y;x_{0} ,y_{0}}  \right)}
\right|}_{x = 0} = 0, \quad {\left. {y^{2\beta} {\frac{{\partial
q_{2} \left( {x,y;x_{0} ,y_{0}} \right)}}{{\partial y}}}}
\right|}_{y = 0} = 0,
\end{equation}

\begin{equation}
\label{eq13} {\left. {x^{2\alpha} {\frac{{\partial q_{3} \left(
{x,y;x_{0} ,y_{0}} \right)}}{{\partial x}}}} \right|}_{x = 0} = 0,
\quad {\left. {q_{3} \left( {x,y;x_{0} ,y_{0}}  \right)}
\right|}_{y = 0} = 0,
\end{equation}

 ${\left. {q_{4} \left( {x,y;x_{0} ,y_{0}}  \right)} \right|}_{x = 0} = 0$ and
${\left. {q_{4} \left( {x,y;x_{0} ,y_{0}}  \right)} \right|}_{y =
0} = 0.$ (1.14)

In the paper [33] using fundamental solution $q_{1} \left(
{x,y;x_{0} ,y_{0} } \right)$ in the domain defined by

\begin{equation}
\label{eq14} \Omega \subset R_{ +} ^{2} = {\left\{ {\left( {x,y}
\right):\,x > 0,\,y > 0} \right\}},
\end{equation}

\noindent the double-layer potential theory for the equation
$\left( {H_{\alpha ,\beta }^{0}}  \right)$ was investigated.

Here, in this publication, we aim at constructing theory of
double-layer potentials corresponding to the next fundamental
solution $q_{2} \left( {x,y;x_{0} ,y_{0}}  \right)$. By using some
properties of one of Appell's hypergeometric functions in two
variables, we prove limiting theorems and derive integral
equations concerning a denseness of double-layer potentials.

\section{Green's formula}

We begin by considering the following identity:

\begin{equation}
\label{eq15} x^{2\alpha} y^{2\beta} {\left[ {uH_{\alpha ,\beta}
^{0} \left( {v} \right) - vH_{\alpha ,\beta} ^{0} \left( {u}
\right)} \right]} = {\frac{{\partial }}{{\partial x}}}{\left[
{x^{2\alpha} y^{2\beta} \left( {v_{x} u - vu_{x}} \right)}
\right]} + {\frac{{\partial}} {{\partial y}}}{\left[ {x^{2\alpha
}y^{2\beta} \left( {v_{y} u - vu_{y}}  \right)} \right]}.
\end{equation}

Integrating both parts of identity (\ref{eq15}) on a domain
$\Omega $ in (\ref{eq14}), and using Green's formula, we find that

\begin{equation}
\label{eq16} {\int\!\!\!\int\limits_{\Omega}  {}} x^{2\alpha}
y^{2\beta} {\left[ {uH_{\alpha ,\beta} ^{0} \left( {v} \right) -
vH_{\alpha ,\beta} ^{0} \left( {u} \right)} \right]}dxdy =
{\int\limits_{S} {}} x^{2\alpha} y^{2\beta }u\left( {v_{x} dy -
v_{y} dx} \right) - x^{2\alpha} y^{2\beta} v\left( {u_{x} dy -
u_{y} dx} \right),
\end{equation}

\noindent where $S = \partial \Omega $ is a boundary of the domain
$\Omega $.

If $u\left( {x,y} \right)$ and $v\left( {x,y} \right)$ are
solutions of the equation$\left( {H_{\alpha ,\beta} ^{0}}
\right)$, we find from (\ref{eq16}) that

\begin{equation}
\label{eq17} {\int\limits_{S} {}} x^{2\alpha} y^{2\beta} \left(
{u{\frac{{\partial v}}{{\partial n}}} - v{\frac{{\partial
u}}{{\partial n}}}} \right)ds = 0,
\end{equation}

\noindent where

\begin{equation}
\label{eq18} {\frac{{\partial}} {{\partial n}}} =
{\frac{{dy}}{{ds}}}{\frac{{\partial }}{{\partial x}}} -
{\frac{{dx}}{{ds}}}{\frac{{\partial}} {{\partial y}}}, \quad
{\frac{{dy}}{{ds}}} = \cos \left( {n,x} \right),
{\frac{{dx}}{{ds}}} = - \cos \left( {n,y} \right),
\end{equation}

 $n$ being the exterior normal to the curve $S$. We also obtain the following
identity:

\begin{equation}
\label{eq19} {\int\!\!\!\int\limits_{\Omega}  {}} x^{2\alpha}
y^{2\beta} {\left[ {u_{x}^{2} + u_{y}^{2} + \lambda ^{2}u^{2}}
\right]}dxdy = {\int\limits_{S} {}} x^{2\alpha} y^{2\beta}
u{\frac{{\partial u}}{{\partial n}}}ds,
\end{equation}

\noindent where $u\left( {x,y} \right)$ is a solution of the
equation $\left( {H_{\alpha ,\beta} ^{0}}  \right)$. The special
case of (\ref{eq17}) when $v = 1$ reduces to the following form:

\begin{equation}
\label{eq20} {\int\limits_{S} {}} x^{2\alpha} y^{2\beta}
{\frac{{\partial u}}{{\partial n}}}ds = 0.
\end{equation}

We note from (\ref{eq20}) that the integral of the normal
derivative of a solution of the equation $\left( {H_{\alpha
,\beta} ^{0}}  \right)$ with a weight $x^{2\alpha} y^{2\beta} $
along the boundary $S$ of the domain $\Omega $ in (\ref{eq14}) is
equal to zero.

\section{A double layer potential $w^{( {2} )}\left(
{x_{0} ,y_{0}}  \right)$}

Let $\Omega $ in (\ref{eq14}) be a domain bounded by intervals
$\left( {0,a} \right)$ and $\left( {0,b} \right)$ of the $x - $
and $y - $ axes, respectively, and a curve $\Gamma $ with the
extremities at points $A( a,0)$ and $B(0,b)$. The parametrical
equations of the curve $\Gamma $ are given by

 $x = x(s)$ and $y = y(s)$ ($s \in [0,l])$,

where $l$ denotes the length of $\Gamma $. We assume the following
properties of the curve $\Gamma $:

(i) The functions $x = x( {s} )$ and $y = y( {s})$ have continuous
derivatives $x'( {s} )$ and $y'( {s} )$ on a segment ${\left[
{0,l} \right]}$, do not vanish simultaneously;

(ii) The second derivatives $x''\left( {s} \right)$ and $y''\left(
{s} \right)$ satisfy to Hoelder condition on ${\left[ {0,l}
\right]}$, where $l$ denotes the length of the curve $\Gamma $;

(iii) In some neighborhoods of points $A\left( {a,0} \right)$ and
$B\left( {0,b} \right)$, the following conditions are satisfied:

\begin{equation} \label{eq001}
 {\left| {{\frac{{dx}}{{ds}}}} \right|} \le cy^{1 + \varepsilon} \left( {s}
\right), \,\,\, {\left| {{\frac{{dy}}{{ds}}}} \right|} \le cx^{1 +
\varepsilon }\left( {s} \right),0 < \varepsilon < 1,\,\,c =
constant
\end{equation}

 $\left( {x,y} \right)$ being the coordinates of a variable point on the curve
$\Gamma $.

Consider the following integral

\begin{equation}
\label{eq21} w^{\left( {2} \right)}\left( {x_{0} ,y_{0}}  \right)
= {\int\limits_{0}^{l} {x^{2\alpha} y^{2\beta}} } \mu _{2} \left(
{s} \right){\frac{{\partial q_{2} \left( {x,y;x_{0} ,y_{0}}
\right)}}{{\partial n}}}ds
\end{equation}

\noindent where the density $\mu _{2} \left( {s} \right) \in
C{\left[ {0,\,l} \right]}$ and $q_{2} $ is given in (\ref{eq3}).
We call the integral (\ref{eq21}) \textit{a double-layer}
\textit{potential} with denseness$\mu _{2} \left( {s} \right)$.

We now investigate some properties of a double-layer potential
$w^{\left( {2} \right)}\left( {x_{0} ,y_{0}}  \right)$ with
denseness $\mu _{2} \left( {s} \right)$.

\textbf{Lemma 1.} \textit{The following formula holds true:}

\begin{equation}
\label{eq22} w_{1}^{\left( {2} \right)} \left( {x_{0} ,y_{0}}
\right) = {\left\{ {{\begin{array}{*{20}c}
 {i\left( {x_{0} ,y_{0}}  \right) - 1} \hfill & {\left( {\left( {x_{0}
,y_{0}}  \right) \in \Omega}  \right)} \hfill \\
 {i\left( {x_{0} ,y_{0}}  \right) - {\frac{{1}}{{2}}}} \hfill & {\left(
{\left( {x_{0} ,y_{0}}  \right) \in \Gamma}  \right)} \hfill \\
 {i\left( {x_{0} ,y_{0}}  \right)} \hfill & {\left( {\left( {x_{0} ,y_{0}}
\right) \notin \bar {\Omega} \,} \right),} \hfill \\
\end{array}}}  \right.}
\end{equation}

\textit{where a domain} $\Omega $\textit{ and the curve} $\Gamma
$\textit{ are described as in this section and} $\bar {\Omega} : =
\Omega \cup \Gamma $;

\[
i\left( {x_{0} ,y_{0}}  \right) = k_{2} \left( {1 - 2\alpha}
\right)x_{0}^{1 - 2\alpha}  {\int\limits_{0}^{b} {}} y^{2\beta}
\left( {x_{0}^{2} + (y - y_{0} )^{2}} \right)^{\alpha - \beta -
1}F\left( {1 - \alpha + \beta ,\beta ;2\beta ;{\frac{{ - 4yy_{0}}}
{{x_{0}^{2} + (y - y_{0} )^{2}}}}} \right)dy.
\]

\textit{Proof.}

\textbf{Case 1.} When $\left( {x_{0} ,y_{0}}  \right) \in \Omega
$, we cut a circle centered at $\left( {x_{0} ,y_{0}}  \right)$
with a small radius $\rho $ off the domain $\Omega $ and denote
the remaining by $\Omega ^{\rho }$ and circuit of the
cut-off-circle by $C_{\rho}  $. The function $q_{2} \left(
{x,y;x_{0} ,y_{0}}  \right)$ in (\ref{eq3}) is a regular solution
of the equation $\left( {H_{\alpha ,\beta} ^{0}}  \right)$ in the
domain $\Omega ^{\rho} $. Using the following derivative formula
of Appell's hypergeometric function ([32], p. 19, (20)):

\begin{equation}
\label{eq23} {\frac{{\partial ^{m + n}F_{2} \left( {a;b_{1} ,b_{2}
;c_{1} ,c_{2} ;x,y} \right)}}{{\partial x^{m}\partial y^{n}}}} =
{\frac{{\left( {a} \right)_{m + n} \left( {b_{1}}  \right)_{m}
\left( {b_{2}}  \right)_{n}}} {{\left( {c_{1} } \right)_{m} \left(
{c_{2}}  \right)_{n}}} }F_{2} \left( {a + m + n;b_{1} + m,b_{2} +
n;c_{1} + m,c_{2} + n;x,y} \right)
\end{equation}

 we have

\[
\begin{array}{l}
 {\frac{{\partial q_{2} \left( {x,y;x_{0} ,y_{0}}  \right)}}{{\partial x}}}
= (1 - 2\alpha )k_{2} \left( {r^{2}} \right)^{\alpha - \beta -
1}x^{ - 2\alpha} x_{0}^{1 - 2\alpha}  F_{2} \left( {1 - \alpha +
\beta ;1 - \alpha
,\beta ;2 - 2\alpha ,2\beta ;\xi ,\eta}  \right) \\
 + 2(\alpha - \beta - 1)k_{2} \left( {r^{2}} \right)^{\alpha - \beta - 2}(x
- x_{0} )x^{1 - 2\alpha} x_{0}^{1 - 2\alpha}  F_{2} \left( {1 -
\alpha +
\beta ;1 - \alpha ,\beta ;2 - 2\alpha ,2\beta ;\xi ,\eta}  \right) \\
 \end{array}
\]

\begin{equation}
\label{eq24}
\begin{array}{l}
 - k_{2} \left( {r^{2}} \right)^{\alpha - \beta - 1}x^{1 - 2\alpha} x_{0}^{1
- 2\alpha}  {\frac{{(1 - \alpha + \beta )(1 - \alpha )}}{{2 -
2\alpha }}}{\frac{{4x_{0}}} {{r^{2}}}}F_{2} \left( {2 - \alpha +
\beta ;2 - \alpha
,\beta ;3 - 2\alpha ,2\beta ;\xi ,\eta}  \right) \\
 - 2k_{2} \left( {r^{2}} \right)^{\alpha - \beta - 2}x^{1 - 2\alpha
}x_{0}^{1 - 2\alpha}  (x - x_{0} ){\left[ {{\frac{{(1 - \alpha +
\beta )(1 - \alpha )}}{{2 - 2\alpha}} }\xi F_{2} \left( {2 -
\alpha + \beta ;2 - \alpha
,\beta ;3 - 2\alpha ,2\beta ;\xi ,\eta}  \right)} \right.} \\
 + {\left. {{\frac{{(1 - \alpha + \beta )\beta}} {{2\beta}} }\eta F_{2}
\left( {2 - \alpha + \beta ;1 - \alpha ,1 + \beta ;2 - 2\alpha ,1
+ 2\beta
;\xi ,\eta}  \right)} \right]} \\
 \end{array}
\end{equation}

By applying the following known contiguous relation (see [32], p.
21):

\begin{equation}
\label{eq25}
\begin{array}{l}
 {\frac{{b_{1}}} {{c_{1}}} }xF_{2} \left( {a + 1;b_{1} + 1,b_{2} ;c_{1} +
1,c_{2} ;x,y} \right) + {\frac{{b_{2}}} {{c_{2}}} }yF_{2} \left(
{a +
1;b_{1} ,b_{2} + 1;c_{1} ,c_{2} + 1;x,y} \right) \\
 = F_{2} \left( {a + 1;b_{1} ,b_{2} ;c_{1} ,c_{2} ;x,y} \right) - F_{2}
\left( {a;b_{1} ,b_{2} ;c_{1} ,c_{2} ;x,y} \right), \\
 \end{array}
\end{equation}

\noindent to (\ref{eq24}), we obtain

\begin{equation}
\label{eq26}
\begin{array}{l}
 {\frac{{\partial q_{2} \left( {x,y;x_{0} ,y_{0}}  \right)}}{{\partial x}}}
= (1 - 2\alpha )k_{2} \left( {r^{2}} \right)^{\alpha - \beta -
1}x^{ - 2\alpha} x_{0}^{1 - 2\alpha}  F_{2} \left( {1 - \alpha +
\beta ;1 - \alpha
,\beta ;2 - 2\alpha ,2\beta ;\xi ,\eta}  \right) \\
 - 2(1 - \alpha + \beta )k_{2} \left( {r^{2}} \right)^{\alpha - \beta -
2}x^{1 - 2\alpha} x_{0}^{2 - 2\alpha}  F_{2} \left( {2 - \alpha +
\beta ;2 -
\alpha ,\beta ;3 - 2\alpha ,2\beta ;\xi ,\eta}  \right) \\
 - 2(1 - \alpha + \beta )k_{2} \left( {r^{2}} \right)^{\alpha - \beta -
2}x^{1 - 2\alpha} x_{0}^{1 - 2\alpha}  (x - x_{0} )F_{2} \left( {2
- \alpha
+ \beta ;1 - \alpha ,\beta ;2 - 2\alpha ,2\beta ;\xi ,\eta}  \right) \\
 \end{array}
\end{equation}

Similarly, we find that

$$
{\frac{{\partial q_{2} \left( {x,y;x_{0} ,y_{0}}
\right)}}{{\partial y}}} = - 2(1 - \alpha + \beta )k_{2} \left(
{r^{2}} \right)^{\alpha - \beta - 2}x^{1 - 2\alpha} x_{0}^{1 -
2\alpha}  y_{0} F_{2} \left( {2 - \alpha + \beta ;1 - \alpha ,1 +
\beta ;2 - 2\alpha ,1 + 2\beta ;\xi ,\eta}  \right)
$$

\begin{equation}
\label{eq27}
 - 2(1 - \alpha + \beta )k_{2} \left( {r^{2}} \right)^{\alpha - \beta -
2}x^{1 - 2\alpha} x_{0}^{1 - 2\alpha}  (y - y_{0} )F_{2} \left( {2
- \alpha + \beta ;1 - \alpha ,\beta ;2 - 2\alpha ,2\beta ;\xi
,\eta}  \right)
\end{equation}

Thus, with the help of (\ref{eq26}) and (\ref{eq27}), it follows
from (\ref{eq3}) and (\ref{eq18}) that

\begin{equation}
\label{eq28}
\begin{array}{l}
 {\frac{{\partial q_{2} \left( {x,y;x_{0} ,y_{0}}  \right)}}{{\partial n}}}
= \\
 = - (1 - \alpha + \beta )k_{2} \left( {r^{2}} \right)^{\alpha - \beta -
1}x^{1 - 2\alpha} x_{0}^{1 - 2\alpha}  F_{2} \left( {2 - \alpha +
\beta ;1 - \alpha ,\beta ;2 - 2\alpha ,2\beta ;\xi ,\eta}
\right){\frac{{\partial
}}{{\partial n}}}{\left[ {\ln r^{2}} \right]} \\
 - 2k_{2} (1 - \alpha + \beta )\left( {r^{2}} \right)^{\alpha - \beta -
2}x^{1 - 2\alpha} x_{0}^{2 - 2\alpha}  F_{2} \left( {2 - \alpha +
\beta ;2 - \alpha ,\beta ;3 - 2\alpha ,2\beta ;\xi ,\eta}
\right){\frac{{dy}}{{ds}}}
\\
 + 2k_{2} (1 - \alpha + \beta )\left( {r^{2}} \right)^{\alpha - \beta -
2}x^{1 - 2\alpha} x_{0}^{1 - 2\alpha}  y_{0} F_{2} \left( {2 -
\alpha + \beta ;1 - \alpha ,1 + \beta ;2 - 2\alpha ,1 + 2\beta
;\xi ,\eta}
\right){\frac{{dx}}{{ds}}} \\
 + (1 - 2\alpha )k_{2} \left( {r^{2}} \right)^{\alpha - \beta - 1}x^{ -
2\alpha} x_{0}^{1 - 2\alpha}  F_{2} \left( {1 - \alpha + \beta ;1
- \alpha
,\beta ;2 - 2\alpha ,2\beta ;\xi ,\eta}  \right){\frac{{dy}}{{ds}}} \\
 \end{array}
\end{equation}

Applying (\ref{eq20}) and considering identity (\ref{eq12}), we
get the following formula:

\begin{equation}
\label{eq29} w_{1}^{\left( {2} \right)} \left( {x_{0} ,y_{0}}
\right) = {\mathop {\lim }\limits_{\rho \to 0}}
{\int\limits_{C_{\rho}}   {}} x^{2\alpha} y^{2\beta
}{\frac{{\partial q_{2} \left( {x,y;x_{0} ,y_{0}}
\right)}}{{\partial n}}}ds + {\int\limits_{0}^{b} {y^{2\beta}
{\left. {{\left[ {x^{2\alpha }{\frac{{\partial q_{2} \left(
{x,y;x_{0} ,y_{0}}  \right)}}{{\partial n}}}} \right]}}
\right|}_{x = 0} ds}} .
\end{equation}

Substituting from (\ref{eq28}) into (\ref{eq29}), we find that

\begin{equation}
\label{eq30}
w_1^{(2)}\left(x_0,y_0\right)=k_2x_0^{1-2\alpha}{\mathop {\lim
}\limits_{\rho \to
0}}\left\{(1-\alpha+\beta)\left[-J_1-2x_0J_2+2y_0J_3\right]+(1-2\alpha)J_4\right\}+J_5,
\end{equation}

\noindent where

$$
J_{1} (x_{0} ,y_{0} ) = {\int\limits_{C_{\rho}}   {}} xy^{2\beta}
\left( {r^{2}} \right)^{\alpha - \beta - 1}F_{2} \left( {2 -
\alpha + \beta ;1 - \alpha ,\beta ;2 - 2\alpha ,2\beta ;\xi ,\eta}
\right){\frac{{\partial }}{{\partial n}}}{\left[ {\ln r^{2}}
\right]}ds,
$$
$$
 J_{2} (x_{0} ,y_{0} ) = {\int\limits_{C_{\rho}} {}} xy^{2\beta}
\left( {r^{2}} \right)^{\alpha - \beta - 2}F_{2} \left( {2 -
\alpha + \beta ;2 - \alpha ,\beta ;3 - 2\alpha ,2\beta ;\xi ,\eta}
\right){\frac{{dy}}{{ds}}}ds,
$$
$$J_{3} (x_{0} ,y_{0} ) = {\int\limits_{C_{\rho}}   {}} xy^{2\beta}
\left( {r^{2}} \right)^{\alpha - \beta - 2}F_{2} \left( {2 -
\alpha + \beta ;1 - \alpha ,1 + \beta ;2 - 2\alpha ,1 + 2\beta
;\xi ,\eta} \right){\frac{{dx}}{{ds}}}ds,
$$
$$J_{4} (x_{0} ,y_{0} ) = \int\limits_{C_{\rho}} y^{2\beta} \left(
{r^{2}} \right)^{\alpha - \beta - 1}F_{2} \left( {1 - \alpha +
\beta ;1 - \alpha ,\beta ;2 - 2\alpha ,2\beta ;\xi ,\eta}
\right){\frac{{dy}}{{ds}}}ds,
$$
$$
J_{5} (x_{0} ,y_{0} ) = {\int\limits_{0}^{b} {y^{2\beta} {\left.
{{\left[ {x^{2\alpha} {\frac{{\partial q_{2} \left( {x,y;x_{0}
,y_{0}} \right)}}{{\partial n}}}} \right]}} \right|}_{x = 0} ds}}.
$$

Now, by introducing the polar coordinates: $x = x_{0} + \rho\cos
\phi $ and $y = y_{0} + \rho \sin \phi $, we get
\begin{equation}
\label{eq31}
\begin{array}{l}
 J_{1} \left( {x_{0} ,y_{0}}  \right) = {\int\limits_{0}^{2\pi}  {(x_{0} +
\rho \cos \phi )(y_{0} + \rho \sin \phi )^{2\beta} \left( {\rho
^{2}}
\right)^{\alpha - \beta - 1}}}  \\
 F_{2} \left( {2 - \alpha + \beta ;1 - \alpha ,\beta ;2 - 2\alpha ,2\beta
;\xi ,\eta}  \right)d\phi \\
 \end{array}
\end{equation}

By using the following known formulas (see [30], p. 253, (26),
[31], p. 113, (4)):

\begin{equation}
\label{eq32} F_{2} \left( {a;b_{1} ,b_{2} ;c_{1} ,c_{2} ;x,y}
\right) = {\sum\limits_{i = 0}^{\infty}  {}} {\frac{{\left( {a}
\right)_{i} \left( {b_{1}}  \right)_{i} \left( {b_{2}}
\right)_{i}}} {{\left( {c_{1}}  \right)_{i} \left( {c_{2}}
\right)_{i} i!}}}x^{i}y^{i}F\left( {a + i,b_{1} + i;c_{1} + i;x}
\right)F\left( {a + i,b_{2} + i;c_{2} + i;y} \right),
\end{equation}
and
\begin{equation}
\label{eq33} F\left( {a,b;c,x} \right) = \left( {1 - x} \right)^{
- b}F\left( {c - a,b;c,{\frac{{x}}{{x - 1}}}} \right),
\end{equation}
we obtain
\begin{equation}
\label{eq34}
\begin{array}{l}
 F_{2} \left( {a;b_{1} ,b_{2} ;c_{1} ,c_{2} ;x,y} \right) = \left( {1 - x}
\right)^{ - b_{1}} \left( {1 - y} \right)^{ - b_{2}}
{\sum\limits_{i = 0}^{\infty}  {}} {\frac{{\left( {a} \right)_{i}
\left( {b_{1}}  \right)_{i} \left( {b_{2}}  \right)_{i}}} {{\left(
{c_{1}}  \right)_{i} \left( {c_{2}} \right)_{i} i!}}}\left(
{{\frac{{x}}{{1 - x}}}} \right)^{i}\left(
{{\frac{{y}}{{1 - y}}}} \right)^{i} \cdot \\
 \cdot F\left( {c_{1} - a,b_{1} + i;c_{1} + i;{\frac{{x}}{{x - 1}}}}
\right)F\left( {c_{2} - a,b_{2} + i;c_{2} + i;{\frac{{y}}{{y -
1}}}}
\right), \\
 \end{array}
\end{equation}
where $F\left( {a,b;c;x} \right)$ is hypergeometric function of
Gauss ([31], p. 69, (2)). Hence we have

\begin{equation}
\label{eq35}
\begin{array}{l}
 F_{2} \left( {2 - \alpha + \beta ;1 - \alpha ,\beta ;2 - 2\alpha ,2\beta
;\xi ,\eta}  \right) \\
 = \left( {\rho ^{2}} \right)^{1 - \alpha + \beta} \left( {\rho ^{2} +
4x_{0}^{2} + 4x_{0} \rho \cos \,\phi}  \right)^{\alpha - 1}\left(
{\rho ^{2}
+ 4y_{0}^{2} + 4y_{0} \rho \sin \,\phi}  \right)^{ - \beta} P_{11} , \\
 \end{array}
\end{equation}
where
\[
\begin{array}{l}
 P_{11} = {\sum\limits_{i = 0}^{\infty}  {}} {\frac{{\left( {2 - \alpha +
\beta}  \right)_{i} \left( {1 - \alpha}  \right)_{i} \left(
{\beta} \right)_{i}}} {{\left( {2 - 2\alpha}  \right)_{i} \left(
{2\beta} \right)_{i} i!}}}\left( {{\frac{{4x_{0}^{2} + 4x_{0} \rho
\cos \,\phi }}{{\rho ^{2} + 4x_{0}^{2} + 4x_{0} \rho \cos \,\phi}}
}} \right)^{i}\left( {{\frac{{4y_{0}^{2} + 4y_{0} \rho \sin
\,\phi}} {{\rho ^{2} + 4y_{0}^{2} +
4y_{0} \rho \sin \,\phi}} }} \right)^{i} \cdot \\
 \cdot F\left( { - \alpha - \beta ,1 - \alpha + i;2 - 2\alpha +
i;{\frac{{4x_{0}^{2} + 4x_{0} \rho \cos \,\phi}} {{\rho ^{2} +
4x_{0}^{2} + 4x_{0} \rho \cos \,\phi}} }} \right)F\left( {\alpha +
\beta - 2,\beta + i;2\beta + i;{\frac{{4y_{0}^{2} + 4y_{0} \rho
\sin \,\phi}} {{\rho ^{2} +
4y_{0}^{2} + 4y_{0} \rho \sin \,\phi}} }} \right). \\
 \end{array}
\]

Using the well-known Gauss's summation formula for $F\left(
{a,b;c;1} \right)$ ([31], p. 112, (46))

\[
F\left( {a,b;c;1} \right) = {\frac{{\Gamma \left( {c}
\right)\Gamma \left( {c - a - b} \right)}}{{\Gamma \left( {c - a}
\right)\Gamma \left( {c - b} \right)}}},c \ne 0, - 1, -
2,...,Re\left( {c - a - b} \right) > 0,
\]
we obtain

\begin{equation}
\label{eq36} {\mathop {\lim} \limits_{\rho \to 0}} P_{11} =
{\frac{{\Gamma (2 - 2\alpha )\Gamma (2\beta )}}{{\Gamma (2 -
\alpha + \beta )\Gamma (\beta )\Gamma (1 - \alpha )}}}.
\end{equation}

Thus, by virtue of the identities (\ref{eq31}), (\ref{eq35}), and
(\ref{eq36}), we get

\begin{equation}
\label{eq37}
 - (1 - \alpha + \beta )k_{2} x_{0}^{1 - 2\alpha}  {\mathop {\lim
}\limits_{\rho \to 0}} J_{1} \left( {x_{0} ,y_{0}}  \right) = - 1.
\end{equation}

Similarly, by considering the corresponding identities and the
fact that

\begin{equation}
\label{eq38} {\mathop {\lim} \limits_{\rho \to 0}} \rho \ln \rho =
0,
\end{equation}
we find that

\begin{equation}
\label{eq39} {\mathop {\lim} \limits_{\rho \to 0}} J_{2} \left(
{x_{0} ,y_{0}}  \right) = {\mathop {\lim} \limits_{\rho \to 0}}
J_{3} \left( {x_{0} ,y_{0}}  \right) = {\mathop {\lim}
\limits_{\rho \to 0}} J_{4} \left( {x_{0} ,y_{0}}  \right) = 0.
\end{equation}

Now we consider the integral $J_{5} \left( {x_{0} ,y_{0}}
\right)$, which, taking into account formula (\ref{eq28}), takes
the form

\begin{equation}
\label{eq40} J_{5} \left( {x_{0} ,y_{0}}  \right) = i(x_{0} ,y_{0}
)
\end{equation}

Hence, by virtue of (\ref{eq37})-(\ref{eq40}), from (\ref{eq30})
at $\left( {x_{0} ,y_{0}} \right) \in \Omega $ follows

\begin{equation}
\label{eq41} w_{1}^{\left( {2} \right)} \left( {x_{0} ,y_{0}}
\right) = i(x_{0} ,y_{0} ) - 1.
\end{equation}

\textbf{Case 2.} When $\left( {x_{0} ,y_{0}}  \right) \in \Gamma
$, we cut a circle $C_{\rho}  $ centered at $\left( {x_{0} ,y_{0}}
\right)$ with a small radius $\rho $ off the domain $\Omega $ and
denote the remaining part of the curve by $\Gamma - \Gamma _{\rho}
$. Let $C_{\rho} ^{'} $ denote a part of the circle $C_{\rho}  $
lying inside the domain $\Omega $. We consider the domain $\Omega
_{\rho}  $which is bounded by a curve $\Gamma - \Gamma _{\rho}  $,
$C_{\rho} ^{'} $ and segments ${\left[ {0,a} \right]}$ and
${\left[ {0,b} \right]}$ along the $x - $ and $y - $axes,
respectively. Then we have

\begin{equation}
\label{eq42}
\begin{array}{l}
 w_{1}^{\left( {2} \right)} \left( {x_{0} ,y_{0}}  \right) \equiv
{\int\limits_{0}^{l} {x^{2\alpha} y^{2\beta}} } {\frac{{\partial
q_{2} \left( {x,y;x_{0} ,y_{0}}  \right)}}{{\partial n}}}ds
 = {\mathop {\lim} \limits_{\rho \to 0}} {\int\limits_{\Gamma - \Gamma
_{\rho}}   {x^{2\alpha} y^{2\beta} {\frac{{\partial q_{2} \left(
{x,y;x_{0}
,y_{0}}  \right)}}{{\partial n}}}ds}}  \\
 \end{array}
\end{equation}

When the point $\left( {x_{0} ,y_{0}}  \right)$ lies outside the
domain $\Omega _{\rho}  $, it is found that, in this domain $q_{2}
\left( {x,y;x_{0} ,y_{0}}  \right)$ is a regular solution of the
equation$\left( {H_{\alpha ,\beta} ^{0}}  \right)$. Therefore, by
virtue of (\ref{eq20}), we have

\begin{equation}
\label{eq43}
\begin{array}{l}
 {\int\limits_{\Gamma - \Gamma _{\rho}}   {x^{2\alpha} y^{2\beta
}{\frac{{\partial q_{2} \left( {x,y;x_{0} ,y_{0}}
\right)}}{{\partial
n}}}ds}}  = \\
 = {\int\limits_{0}^{b} {y^{2\beta}} } {\left. {{\left[ {x^{2\alpha
}{\frac{{\partial q_{2} \left( {x,y;x_{0} ,y_{0}}
\right)}}{{\partial x}}}} \right]}} \right|}_{x = 0} dy +
{\int\limits_{C_{\rho} ^{'}}  {}} x^{2\alpha }y^{2\beta}
{\frac{{\partial}} {{\partial n}}}{\left\{ {q_{2} \left(
{x,y;x_{0} ,y_{0}}  \right)} \right\}}ds. \\
 \end{array}
\end{equation}

Substituting from (\ref{eq43}) into (\ref{eq42}), we get

\begin{equation}
\label{eq44} w_{1}^{\left( {2} \right)} \left( {x_{0} ,y_{0}}
\right) = i(x_{0} ,y_{0} ) + {\mathop {\lim} \limits_{\rho \to 0}}
{\int\limits_{C_{\rho} ^{'}}  {} }x^{2\alpha} y^{2\beta}
{\frac{{\partial q_{2} \left( {x,y;x_{0} ,y_{0}}
\right)}}{{\partial n}}}ds.
\end{equation}

Now, again by introducing the polar coordinates in the second
summand and calculating the limit at $\rho \to 0$, we obtain

\[
w_{1}^{\left( {2} \right)} \left( {x_{0} ,y_{0}}  \right) =
i(x_{0} ,y_{0} ) - {\frac{{1}}{{2}}}.
\]

\textbf{Case 3.} When $\left( {x_{0} ,y_{0}}  \right) \notin \bar
{\Omega }$, it is noted that the function $q_{2} \left( {x,y;x_{0}
,y_{0}}  \right)$ is a regular solution of the equation $\left(
{H_{\alpha ,\beta} ^{0}} \right)$. Hence, in view of the formula
(\ref{eq20}), we have

\[
\begin{array}{l}
 w_{1}^{\left( {2} \right)} \left( {x_{0} ,y_{0}}  \right) \equiv
{\int\limits_{0}^{l} {x^{2\alpha} y^{2\beta}} } {\frac{{\partial
}}{{\partial n}}}{\left\{ {q_{2} \left( {x,y;x_{0} ,y_{0}}
\right)} \right\}}ds
 = {\int\limits_{0}^{b} {y^{2\beta}} } {\left. {{\left[ {x^{2\alpha
}{\frac{{\partial q_{2} \left( {x,y;x_{0} ,y_{0}}
\right)}}{{\partial x}}}}
\right]}} \right|}_{x = 0} dy = i(x_{0} ,y_{0} ). \\
 \end{array}
\]

The proof of Lemma 1 is thus completed.

\textbf{Lemma 2.} \textit{The following formula holds true:}
\begin{equation}
\label{eq45} w_{1}^{\left( {2} \right)} \left( {x_{0} ,0} \right)
= {\left\{ {{\begin{array}{*{20}c}
 {i\left( {x_{0} ,0} \right) - 1} \hfill & {\left( {x_{0} \in \left( {0,a}
\right)} \right)} \hfill \\
 {i\left( {x_{0} ,0} \right) - {\frac{{1}}{{2}}}} \hfill & {\left( {x_{0} =
0\,\,{\rm o}{\rm r}\,\,x_{0} = a} \right)} \hfill \\
 {i\left( {x_{0} ,0} \right)} \hfill & {\left( {a < x_{0}}  \right)} \hfill
\\
\end{array}}}  \right.}
\end{equation}
where
$$
i\left( {x_{0} ,0} \right) = {\frac{{1 - 2\alpha}} {{1 + 2\beta}}
}k_{2} b^{2\beta + 1}x_{0}^{1 - 2\alpha}  \left( {x_{0}^{2} +
b^{2}} \right)^{ - 1 + \alpha - \beta} F\left( {1,\beta +
{\frac{{1}}{{2}}};\beta +
{\frac{{3}}{{2}}};{\frac{{b^{2}}}{{x_{0}^{2} + b^{2}}}}} \right).
$$

\textit{Proof.} For considering the first case when $x_{0} \in
\left( {0,a} \right)$, we introduce a straight line $y = h$ for a
sufficiently small positive real number $h$ and consider a domain
$\Omega _{h} $  which is the part of the domain $\Omega $ lying
above the straight line $y = h$. Applying the formula
(\ref{eq20}), we obtain

\begin{equation}
\label{eq46} w_{1}^{\left( {2} \right)} \left( {x_{0} ,0} \right)
= {\int\limits_{0}^{b} {}} {\left. {x^{2\alpha} y^{2\beta}
{\frac{{\partial q_{2} \left( {x,y;x_{0} ,0} \right)}}{{\partial
x}}}} \right|}_{x = 0} dy + {\mathop {\lim }\limits_{h \to 0}}
{\int\limits_{0}^{x_{1}}  {}} {\left. {x^{2\alpha }y^{2\beta}
{\frac{{\partial q_{2} \left( {x,y;x_{0} ,0} \right)}}{{\partial
y}}}} \right|}_{y = h} dx,
\end{equation}
where $x_{1} \left( {\varepsilon}  \right)$ is an abscissa of a
point at which the straight line $y = h$ intersects the curve
$\Gamma $. It follows from (\ref{eq40}), (\ref{eq27}) and
(\ref{eq46}) that

\begin{equation}
\label{eq47}
\begin{array}{l}
 w_{1}^{\left( {2} \right)} \left( {x_{0} ,0} \right) = i(x_{0} ,0) - \\
 - 2\left( {1 - \alpha + \beta}  \right)k_{2} x_{0}^{1 - 2\alpha}  {\mathop
{\lim} \limits_{h \to 0}} h^{1 + 2\beta} {\int\limits_{0}^{x_{1}}
{} }x{\frac{{F\left( {2 - \alpha + \beta ,1 - \alpha ;2 - 2\alpha
,{\frac{{ - 4xx_{0}}} {{\left( {x - x_{0}}  \right)^{2} +
h^{2}}}}} \right)}}{{{\left[ {\left( {x - x_{0}}  \right)^{2} +
h^{2}} \right]}^{2 - \alpha + \beta
}}}}dx. \\
 \end{array}
\end{equation}

Now, by using the hypergeometric transformation formula
(\ref{eq33}) inside the integrand (\ref{eq47}), we have

\begin{equation}
\label{eq48}
\begin{array}{l}
 w_{1}^{\left( {2} \right)} \left( {x_{0} ,0} \right) = i(x_{0} ,0) - \\
 2\left( {1 - \alpha + \beta}  \right)k_{2} x_{0}^{1 - 2\alpha}  {\mathop
{\lim} \limits_{h \to 0}} h^{1 + 2\beta} {\int\limits_{0}^{x_{1}}
{} }x{\frac{{F\left( { - \alpha - \beta ,1 - \alpha ;2 - 2\alpha
,{\frac{{4xx_{0}}} {{\left( {x + x_{0}}  \right)^{2} + h^{2}}}}}
\right)}}{{{\left[ {\left( {x - x_{0}}  \right)^{2} + h^{2}}
\right]}^{1 + \beta} {\left[ {\left( {x + x_{0}}  \right)^{2} +
h^{2}} \right]}^{1 -
\alpha}} }}dx, \\
 \end{array}
\end{equation}
which, upon setting $x = x_{0} + ht$ inside the integrand, yields

\begin{equation}
\label{eq49}
\begin{array}{l}
 w_{1}^{\left( {2} \right)} \left( {x_{0} ,0} \right) = i(x_{0} ,0) - \\
 - 2\left( {1 - \alpha + \beta}  \right)k_{2} x_{0}^{1 - 2\alpha}  {\mathop
{\lim} \limits_{h \to 0}} {\int\limits_{l_{1}} ^{l_{2}}  {}}
\left( {x_{0} + ht} \right){\frac{{F\left( { - \alpha - \beta ,1 -
\alpha ;2 - 2\alpha ,{\frac{{4x_{0} \left( {x_{0} + ht}
\right)}}{{\left( {2x_{0} + ht} \right)^{2} + h^{2}}}}}
\right)}}{{\left( {1 + t^{2}} \right)^{\beta + 1}{\left[ {\left(
{2x_{0} + ht} \right)^{2} + h^{2}} \right]}^{1 - \alpha
}}}}dt, \\
 \end{array}
\end{equation}
where

$$
l_{1} = - {\frac{{x_{0}}} {{h}}}, \quad l_{2} = {\frac{{x_{1} -
x_{0}}} {{h}}}.
$$

Considering

$$
{\mathop {\lim} \limits_{h \to 0}} F\left( { - \alpha - \beta ,1 -
\alpha ;2 - 2\alpha ,{\frac{{4x_{0} \left( {x_{0} + ht}
\right)}}{{\left( {2x_{0} + ht} \right)^{2} + h^{2}}}}} \right) =
$$
$$
F\left( { - \alpha - \beta ,1 - \alpha ;2 - 2\alpha ,1} \right) =
{\frac{{\Gamma \left( {2 - 2\alpha} \right)\Gamma \left( {1 +
\beta}  \right)}}{{\Gamma \left( {2 - \alpha + \beta}
\right)\Gamma \left( {1 - \alpha}  \right)}}},
$$
and

$$
{\int\limits_{ - \infty} ^{ + \infty}  {}} {\frac{{dt}}{{\left( {1
+ t^{2}} \right)^{\beta + 1}}}} = {\frac{{\pi \Gamma \left(
{2\beta} \right)}}{{2^{2\beta - 1}\beta \Gamma ^{2}\left( {\beta}
\right)}}},
$$
we find from (\ref{eq49}) that

\begin{equation}
\label{eq50} w_{1}^{\left( {2} \right)} \left( {x_{0} ,0} \right)
= i(x_{0} ,0) - 1.
\end{equation}

The other three cases when $x_{0} = 0$, $x_{0} = a$ and$x_{0} > a$
can be proved by using arguments similar to those detailed above
in the first case.

This evidently completes our proof of Lemma 2.

\textbf{Lemma 3.} \textit{The following identities are fair}

\begin{equation}
\label{eq51} w_{1}^{\left( {2} \right)} \left( {0,y_{0}}  \right)
= {\left\{ {{\begin{array}{*{20}c}
 { - 1} \hfill & {\left( {y_{0} \in \left( {0,b} \right)} \right)} \hfill \\
 { - {\frac{{1}}{{2}}}} \hfill & {\left( {y_{0} = 0\,\,{\rm o}{\rm
r}\,\,y_{0} = b} \right)} \hfill \\
 {0} \hfill & {\left( {b < y_{0}}  \right)} \hfill \\
\end{array}}}  \right.}
\end{equation}

\textit{Proof.} The proof of Lemma 3 would run parallel to that of
Lemma 2.

\textbf{Theorem 1.} \textit{For any points} $\left( {x,y}
\right)$\textit{ and} $\left( {x_{0} ,y_{0}}  \right) \in R_{ +}
^{2} $\textit{ and} $x \ne x_{0} $\textit{ and} $y \ne y_{0}
,$\textit{ the following inequality holds true:}

\begin{equation}
\label{eq52} {\left| {q_{2} \left( {x,y;x_{0} ,y_{0}}  \right)}
\right|} \le {\frac{{\Gamma \left( {1 - \alpha}  \right)\Gamma
\left( {\beta} \right)}}{{\pi \Gamma (1 - \alpha + \beta
)}}}{\frac{{4^{\beta - \alpha }x^{1 - 2\alpha} x_{0}^{1 -
2\alpha}} } {{\left( {r_{1}^{2}}  \right)^{1 - \alpha} \left(
{r_{2}^{2}}  \right)^{\beta}} }}F{\left[ {1 - \alpha ,\beta ;1 -
\alpha + \beta ;\left( {1 - {\frac{{r^{2}}}{{r_{1}^{2}}} }}
\right)\left( {1 - {\frac{{r^{2}}}{{r_{2}^{2}}} }} \right)}
\right]},
\end{equation}
where $\alpha $\textit{and} $\beta $\textit{ are real parameters
with} $\left( {0 < \alpha ,\beta < {\frac{{1}}{{2}}}}
\right)$\textit{ as in the equation} $\left( {H_{\alpha ,\beta}
^{\lambda}}   \right)$ (\textit{with} $\lambda = 0)$\textit{, and}
$r,r_{1} $\textit{ and} $r_{2} $\textit{ are as in} (\ref{eq10})$.
$

\textit{Proof.} It follows from (\ref{eq34}) that

\begin{equation}
\label{eq53}
\begin{array}{l}
 q_{2} \left( {x,y;x_{0} ,y_{0}}  \right) = k_{2} x^{1 - 2\alpha} x_{0}^{1 -
2\alpha}  \left( {r_{1}^{2}}  \right)^{\alpha - 1}\left(
{r_{2}^{2}} \right)^{ - \beta} {\sum\limits_{i = 0}^{\infty}  {}}
{\frac{{\left( {1 - \alpha + \beta}  \right)_{i} \left( {1 -
\alpha}  \right)_{i} \left( {\beta } \right)_{i}}} {{\left( {2 -
2\alpha}  \right)_{i} \left( {2\beta} \right)_{i} i!}}}\left( {1 -
{\frac{{r^{2}}}{{r_{1}^{2}}} }}
\right)^{i}\left( {1 - {\frac{{r^{2}}}{{r_{2}^{2}}} }} \right)^{i}\times \\
 \times F\left( {1 - \alpha - \beta ,1 - \alpha + i;2 - 2\alpha + i;1 -
{\frac{{r^{2}}}{{r_{1}^{2}}} }} \right)F\left( {\alpha + \beta -
1,\beta +
i;2\beta + i;1 - {\frac{{r^{2}}}{{r_{2}^{2}}} }} \right), \\
 \end{array}
\end{equation}

Now, in view of the following inequalities:

$$
F\left( {1 - \alpha - \beta ,1 - \alpha + i;2 - 2\alpha + i;1 -
{\frac{{r^{2}}}{{r_{1}^{2}}} }} \right) \le {\frac{{(2 - 2\alpha
)_{i} \Gamma (2 - 2\alpha )\Gamma (\beta )}}{{(1 - \alpha + \beta
)_{i} \Gamma (1 - \alpha + \beta )\Gamma (1 - \alpha )}}}
$$
and
$$ F\left( {\alpha + \beta - 1,\beta + i;2\beta + i;1 -
{\frac{{r^{2}}}{{r_{2}^{2}}} }} \right) \le {\frac{{(2\beta )_{i}
\Gamma (2\beta )\Gamma (1 - \alpha )}}{{(1 - \alpha + \beta )_{i}
\Gamma (1 - \alpha + \beta )\Gamma (\beta )}}},
$$
we find from (\ref{eq53}) that the inequality (\ref{eq52}) holds
true. Hence Theorem 1 is proved.

By virtue of the following known formula ([31], p. 117, (12)):
$$
F\left( {a,b;a + b;z} \right) = - {\frac{{\Gamma \left( {a + b}
\right)}}{{\Gamma \left( {a} \right)\Gamma \left( {b}
\right)}}}F\left( {a,b;1;1 - z} \right)\ln \left( {1 - z} \right)
+
$$
$$
 + {\frac{{\Gamma \left( {a + b} \right)}}{{\Gamma ^{2}\left( {a}
\right)\Gamma ^{2}\left( {b} \right)}}}{\sum\limits_{j =
0}^{\infty}  {} }{\frac{{\Gamma \left( {a + j} \right)\Gamma
\left( {b + j} \right)}}{{\left( {j!} \right)^{2}}}}{\left[ {2\psi
\left( {1 + j} \right) - \psi \left( {a + j} \right) - \psi \left(
{b + j} \right)} \right]}\left( {1 - z} \right)^{j},
$$
$$
\left( { - \pi < \arg \,\left( {1 - z} \right) < \pi ,\,\,a,b \ne
0, - 1, - 2,...} \right),
$$
we observe from (\ref{eq52}) that $q_{2} \left( {x,y;x_{0}
,y_{0}}  \right)$ has a logarithmic singularity at $r = 0$.

\textbf{Theorem 2.} \textit{If the curve} $\Gamma $\textit{
satisfies to conditions (3.1) the inequality takes place}

$$
{\int\limits_{\Gamma}  {}} x^{2\alpha} y^{2\beta} {\left|
{{\frac{{\partial q_{2} \left( {x,y;x_{0} ,y_{0}}
\right)}}{{\partial n}}}} \right|}ds \le C_{1} ,
$$
where $C_{1} $\textit{ is a constant.}

\textit{Proof.} Theorem 2 follows by suitably applying Lemmas 1 to
3.

\textbf{Theorem 3.} \textit{The following limiting formulas hold
true for a double-layer potential (\ref{eq21}):}

\begin{equation}
\label{eq54} w_{i}^{\left( {2} \right)} \left( {t} \right) = -
{\frac{{1}}{{2}}}\mu _{2} \left( {t} \right) +
{\int\limits_{0}^{l} {}} \mu _{2} \left( {s} \right)K_{2} \left(
{s,t} \right)ds
\end{equation}
and
\begin{equation} \label{eq55} w_{e}^{\left( {2} \right)} \left(
{t} \right) = {\frac{{1}}{{2}}}\mu _{2} \left( {t} \right) +
{\int\limits_{0}^{l} {}} \mu _{2} \left( {s} \right)K_{2} \left(
{s,t} \right)ds,
\end{equation}
where, as usual, $\mu _{2} \left( {t} \right) \in С{\left[
{0,\,\,l} \right]},$
$$
K_{2} \left( {s,t} \right) = [x(s)]^{2\alpha} [y(s)]^{2\beta
}{\frac{{\partial}} {{\partial n}}}{\left\{ {q_{2} {\left[
{x\left( {s} \right),y\left( {s} \right);x_{0} \left( {t}
\right),y_{0} \left( {t} \right)} \right]}} \right\}}
$$

$$
\left( {\left( {x\left( {s} \right),y\left( {s} \right)} \right)
\in \Gamma ;\left( {x_{0} \left( {t} \right),y_{0} \left( {t}
\right)} \right) \in \Gamma}  \right),
$$
$w_{i}^{\left( {2} \right)} \left( {t} \right)$ and $w_{e}^{\left(
{2} \right)} \left( {t} \right)$ are limiting values of the
double-layer potential (\ref{eq21}) at

$$
\left( {x_{0} \left( {t} \right),y_{0} \left( {t} \right)} \right)
\to \Gamma
$$
from the inside and the outside, respectively.

\textit{Proof.} We find from Lemma 1, in conjunction with Theorems
1 and 2, that each of the limiting formulas asserted by Theorem 3
holds true.

\textbf{References}

[1] A. Altin, Solutions of type $r^{m}$ for a class of singular
equations. \textit{Internat. J. Math. Sci.,} \textbf{5}(1982),
613--619.

[2] A. Altin, Some expansion formulas for a class of singular
partial differential equations. \textit{Proc. Amer. Math. Soc.},
\textbf{85}(1982), 42--46.

[3] A. Altin and E.C.Young, Some properties of solutions of a
class of singular partial differential equations. \textit{Bull.
Ins. Math. Acad. Sinica}, \textbf{11}(1983), 81--87.

[4] P. Appell and J. Kampe de Feriet, \textit{Fonctions
Hypergeometriques et Hyperspheriques; Polynomes d'Hermite,}
Gauthier - Villars. Paris, 1926.

[5] J.L.Burchnall and T.W.Chaundy, Expansions of Appell's double
hypergeometric functions. // \textit{Quart. J. Math. Oxford Ser.}
11(1940), 249-270.

[6] A.Erdelyi, Singularities of generalized axially symmetric
potentials. \textit{Comm. Pure Appl. Math.}, 2(1956), 403-414.

[7] A.Erdelyi, An application of fractional integrals. \textit{J.
Analyse. Math}., 14(1965), 113-126.

[8] A.Erdelyi, W.Magnus, F.Oberhettinger and F.G.Tricomi,
\textit{Higher Transcendental Functions,} Vol.I, McGraw-Hill Book
Company, New York, Toronto and London,1953; Russian edition,
Izdat. Nauka, Moscow, 1973.

[9] A.J.Fryant, Growth and complete sequences of generalized
bi-axially symmetric potentials. \textit{J. Differential
Equations}, \textbf{31}(1979), 155--164.

[10] R.P.Gilbert, On the singularities of generalized axially
symmetric potentials. \textit{Arch. Rational Mech. Anal.}, 6
(1960), 171-176.

[11] R.P.Gilbert, Some properties of generalized axially symmetric
potentials. \textit{Amer. J. Math.}, 84(1962), 475-484.

[12] R.P.Gilbert, ``Bergman's'' integral operator method in
generalized axially symmetric potential theory. \textit{J. Math.
Phys}., 5(1964), 983-987.

[13] R.P.Gilbert, On the location of singularities of a class of
elliptic partial differential equations in four variables.
\textit{Canad. J. Math}., 17(1965), 676-686.

[14] R.P.Gilbert, On the analytic properties of solutions to a
generalized axially symmetric Schroodinger equations. \textit{J.
Differential equations}, 3(1967), 59-77.

[15] R.P.Gilbert, An investigation of the analytic properties of
solutions to the generalized axially symmetric, reduced wave
equation in $n + 1$ variables, with an application to the theory
of potential scattering. \textit{SIAM J. Appl. Math.} 16 (1968),
30-50.

[16] R.P.Gilbert,\textit{ Function Theoretic Methods in Partial
Differential Equations.} New York, London: Academic Press.

[17] R.P.Gilbert and H.Howard, On solutions of the generalized
axially symmetric wave equation represented by Bergman operators,
\textit{Proc. London Math. Soc}., \textbf{15} (1965), 346-360.

[18] N.M.Gunter, \textit{Potential Theory and Its Applications to
Basic Problems of Mathematical Physics} (Translated from the
Russian edition by J.R.Schulenberger), Frederick Ungar Publishing
Company, New York, 1967.

[19] A.Hasanov, Fundamental solutions of generalized bi-axially
symmetric Helmholtz equation. \textit{Complex Variables and
Elliptic Equations}, \textbf{52}(2007), 673--683.

[20] P.Henrici, Zur Funktionentheorie der Wellengleichung,
\textit{Comment. Math. Helv}., 27(1953), 235-293.

[21] P.Henrici, On the domain of regularity of generalized axially
symmetric potentials. \textit{Proc. Amer. Math. Soc}., 8(1957),
29-31.

[22] P.Henrici, Complete systems of solutions for a class of
singular elliptic Partial Differential Equations. \textit{Boundary
Value Problems in differential equations, pp.19-34,University of
Wisconsin Press, Madison}, 1960.

[23] A.Huber, On the uniqueness of generalized axisymmetric
potentials. \textit{Ann. Math}., 60(1954), 351-358.

[24] D.Kumar, Approximation of growth numbers generalized
bi-axially symmetric potentials. \textit{Fasciculi Mathematics},
35(2005), 51--60.

[25] C.Y.Lo, Boundary value problems of generalized axially
symmetric Helmholtz equations. \textit{Portugaliae Mathematica}.
36(1977), 279-289.

[26] O.I.Marichev, An integral representation of solutions of the
generalized biaxiallu symmetric Helmholtz equation and formulas
its inversion (Russian). \textit{Differencial'nye Uravnenija},
Minsk, \textbf{14}(1978), 1824--1831.

[27] P.A.McCoy, Polynomial approximation and growth of generalized
axisymmetric potentials. \textit{Canadian Journal of Mathematics},
\textbf{31}(1979), 49--59.

[28] P.A.McCoy, Best $L^{p}$-approximation of Generalized
bi-axisymmetric Potentials. \textit{Proceedings of the American
Mathematical Society}, \textbf{79}(1980), 435-440.

[29] C.Miranda, \textit{Partial Differential Equations of Elliptic
Type. Second Revised Edition} (Translated from the Italian edition
by Z.C.Motteler), Ergebnisse der Mathematik und ihrer
Grenzgebiete, Band 2, Springer-Verlag, Berlin, Heidelberg and New
York, 1970.

[30] P.P.Niu and X.B.Lo, Some notes on solvability of LPDO.
\textit{Journal of Mathematical Research and Expositions},
\textbf{3}(1983), 127--129.

[31] K.B.Ranger, On the construction of some integral operators
for generalized axially symmetric harmonic and stream
functions.\textit{ J. Math. Mech}., 14(1965), 383-401.

[32] J.M.Rassias and A.Hasanov, Fundamental Solutions of Two
Degenerated Elliptic Equations and Solutions of Boundary Value
Problems in Infinite Area. \textit{International Journal of
Applied Mathematics \& Statistics}, \textbf{8}(2007), 87-95.

[33] H.M.Srivastava, A.Hasanov and J.Choi, 2015. Double-Layer
Potentials for a Generalized Bi-Axially Symmetric Helmholtz
Equation. \textit{Sohag J.Math}. 2, No.1(2015),1-10.

[34] H.M.Srivastava and Karlsson, \textit{Multipl. Gaussian
Hypergeometric Series}, Halsted Press (Ellis Horwood Limited,
Chicherster), John Wiley and Sons, New York,Chichester,Brisbane
and Toronto,1985.

[35] R.J.Weinacht, Some properties of generalized axially
symmetric Helmholtz potentials. \textit{SIAM J. Math. Anal}.
5(1974), 147-152.

[36] A.Weinstein, Discontinuous integrals and generalized
potential theory. \textit{Trans. Amer. Math. Soc}., 63(1948),
342-354.

[37] A.Weinstein, Generalized axially symmetric potentials theory.
\textit{Bull. Amer. Math. Soc}., 59(1953), 20-38.

\end{document}